\documentclass[11pt]{amsart}

\usepackage{amsmath}
\usepackage{amssymb}

\newtheorem{theorem}{Theorem}[section]
\newtheorem{claim}[theorem]{Claim}

\theoremstyle{definition}
\newtheorem{definition}[theorem]{Definition}

\theoremstyle{remark}
\newtheorem{remark}[theorem]{Remark}

\newcount\skewfactor
\def\mathunderaccent#1#2 {\let\theaccent#1\skewfactor#2
\mathpalette\putaccentunder}
\def\putaccentunder#1#2{\oalign{$#1#2$\crcr\hidewidth
\vbox to.2ex{\hbox{$#1\skew\skewfactor\theaccent{}$}\vss}\hidewidth}}
\def\name{\mathunderaccent\tilde-3 }

\newcommand{\pa}{\forall}

\newcommand{\ps}{\subseteq}

\newcommand{\lk}{\langle}
\newcommand{\rk}{\rangle}

\newcommand{\wto}{\Rightarrow}
\newcommand{\uph}{\upharpoonright}

\newcommand{\bbP}{{\mathbb P}}

\newcommand{\sF}{{\mathcal{F}}}

\newcommand{\sP}{{\mathcal{P}}}

\newcommand{\cf}{{\rm cf}}

\newcommand{\st}{{such that}}
\newcommand{\seq}{{sequence}}
\newcommand{\cont}{{continuous}}
\newcommand{\incr}{{increasing}}

\newcommand{\Wlog}{{without loss of generality}}

\title{Weak Diamond}
\author{Saharon Shelah}
\address{Institute of Mathematics
 The Hebrew University of Jerusalem
 Jerusalem 91904, Israel
 and  Department of Mathematics
 Rutgers University
 New Brunswick, NJ 08854, USA}
\email{shelah@math.huji.ac.il}
\urladdr{http://www.math.rutgers.edu/\char`\~shelah}
\thanks{
Research supported by the United States-Israel Binational Science Foundation
Publication 755}

\subjclass{}
\keywords{Set theory, Normal ideals, Weak diamond, precipituous filters,
semi saturated filters }

\begin{document}

\begin{abstract}
Under some cardinal arithmetic assumptions, we prove that every
stationary of $\lambda$ of a right cofinality has the weak
diamond.  This is a strong negation of uniformization. We then
deal with a  weaker version of the weak diamond- colouring
restrictions. We then deal with semi- saturated (normal) filters.
\end{abstract}

\maketitle

\newpage

\underline{Annotated Content}

\S 1.  {\bf{Weak Diamond: sufficient condition}}

[We prove that if $\lambda = 2^{\mu} = \lambda^{< \lambda}$ is
weakly inaccessible,
$$
\Theta= \{\theta:\theta = \cf (\theta) < \lambda
 \quad \hbox{and} \quad \alpha < \lambda \wto
|\alpha|^{\lk {\rm tr},\theta \rk} <\lambda \quad \hbox{and}
\quad S \ps \{\delta<\lambda: \cf (\delta) \in \Theta\}
$$
is stationary then is has weak diamond. We can omit or weaken the
demand $\lambda = \lambda^{< \lambda}$ if we restrict the colouring
${\bf F}$ \st\ for $\eta \in {}^\delta\delta, {\bf F} (\eta)$
depend only on $\eta \uph C_\delta$ where $C_\delta \ps \delta,
\lambda=\lambda^{|C_\delta|}]$.
${}$

\S 2. {\bf{On versions of precipitousness}}

[We show that for successor $\lambda > \beth_\omega$, the
club filter on $\lambda$ is not semi saturated (even
every normal filter concentrating on any $S \in I[\lambda]$ of
cofinality from a large family). Woodin has proved $D_{\omega_2}
+S^2_0$ consistently is semi saturated].
\newpage

\section{Weak Diamond: sufficient condition}
On the weak diamond see \cite{DvSh:65}, \cite[Appendix \S 1]{Sh:f},
\cite{Sh:208}, \cite{Sh:638}; there will be subsequent work on the
middle diamond.

\begin{definition}
\label{x.0}
For regular uncountable $\lambda$,
\begin{enumerate}
\item We say $S \ps \lambda$ is small if
it is ${\bf F}$-small for some function ${\bf F}$ from
${}^{\lambda \ge}\lambda$ to $\{0,1\}$, which means

$(*)_{{\bf F},S}$ for every $\bar{c} \in {}^S 2$ there is $\eta
\in {}^\lambda\lambda$ \st\ $\{ \lambda \in S: {\bf F} (\eta
\uph \delta) = c_\delta \}$ not is stationary.
\item Let $D^{wd}_\lambda = \{A \ps \lambda: \lambda \setminus
A$ is small \}, it is a normal deal (the weak diamond ideal).
\end{enumerate}
\end{definition}

\begin{claim}
\label{x.1}
Assume
\begin{enumerate}
\item[(a)] $\lambda = \lambda^{<\lambda} = 2^\mu$
\item[(b)] $\Theta = \{\theta:\theta= \cf (\theta)$ and for every
$\alpha < \lambda$, we have $|\alpha|^{<\theta>} < \lambda$ or just
$|\alpha|^{< {\rm tr}, \theta>} < \lambda \}$
(see below; so if $\lambda > \beth_\omega$ every large enough
$\theta < \beth_\omega$ is in $\Theta$)
\item[(c)] $\theta \in \Theta$ and $S \ps \{ \delta <
\lambda : \cf(\delta) = \theta, \mu^\omega$ divides $\delta\}$
is stationary.
\end{enumerate}
\underline{Then} $S$ is not in the ideal $D^{wd}_\lambda$ of small
subsets of $\lambda$.
\end{claim}

\begin{definition}
\label{x.2}
\begin{enumerate}
\item Let $\chi^{\lk \theta \rk} = {\rm Min} \{ | \sP| : \sP
\ps [\chi]^\theta$
and every $A \in [\chi]^\theta$ is included in the
union of $< \theta$ members of $\sP \}$.
\item  $\chi^{\lk {\rm tr}, \theta \rk} = \sup
\{ |{\rm lim}_\theta (t) | : t$ is a tree with $\le \chi$
modes and $\theta$ levels \}
\end{enumerate}
\end{definition}

\begin{remark}
\label{x.3r}
\begin{enumerate}
\item On $\chi^{\lk {\rm tr},\theta \rk}$ see \cite{Sh:589},
on $\chi^{\lk \theta\rk}$ see there and in \cite{Sh:460} but no real
knowledge is assumed here.
\item The interesting case of \ref{x.1} is $\lambda$ (weakly)
inaccessible; for $\lambda$ successor we know more; but in later
results even if $2^\mu$ is successor we say on it new things.
\end{enumerate}
\end{remark}
\medskip

\begin{proof}
Let $\bf{F}$ be a function from $\bigcup\limits_{\delta \in S}
{}^{\delta}2$ to $\{0,1\}$, i.e. ${\bf F}$ is a colouring, and we shall
find
$f \in {}^{S}2$ as required for it.

Let $ \{ \nu_i : i < \lambda \}$ list
$\bigcup\limits_{\alpha<\lambda} {}^{\alpha}2$ such that
$$
\alpha < {\rm lg} (\nu_i) \wto \nu_i \uph \alpha \in
\{\nu_j : j < i\}.
$$
For $\delta \in S$  let $\sP_\delta = \{ \eta
\in {}^\delta2$ : \quad $(\pa \alpha < \delta)
(\eta \uph \alpha \in \{\nu_i:i < \delta \})\}.$

Clearly $\delta \in S \wto |\sP_\delta| \le |\delta|^{<{\rm tr},
\theta>} <\lambda$ by assumption (c).
For each $\eta \in \sP_\delta$ we define $h_\eta \in {}^\mu2$
by: $h_\eta (\varepsilon)= {\bf F} (g_{\eta,\varepsilon})$
where for $\varepsilon < \mu$, we let $g_{\eta, \varepsilon}
\in {}^\delta2$ be defined by $g_{\eta,\varepsilon} (\alpha)
= \eta(\mu \alpha+\varepsilon)$ for $\alpha<\delta$, recalling that
$\mu^\omega$ divides $\delta$ as $\delta \in S$.
So $\{h_\eta :\eta \in
\sP_\delta \}$ is a subset of ${}^\mu2$ of cardinality
$\le |\sP_\delta | < \lambda = 2^\mu$ hence we can choose
$g^*_\delta \in {}^\mu2 \setminus \{g_\eta :
\eta \in \sP_\delta\}.$ For $\varepsilon < \mu$ let $f_\varepsilon
\in {}^S2$ be $f_\varepsilon (\delta) = 1 - g^*_\delta
(\varepsilon).$ If for some $\varepsilon <\mu$ the
function $f_\varepsilon$ serve as a weak diamond
\seq\  for ${\bf F}$, we are done so assume that (for each
$\varepsilon <\mu)$ there are
$\eta_\varepsilon$ and $E_\varepsilon$ \st:
\begin{enumerate}
\item[(a)] $E_\varepsilon$ is a club of $\lambda$.
\item[(b)] $\eta_\varepsilon \in {}^\lambda2.$
\item[(c)] if $\delta \in E_\varepsilon \cap S$
then ${\bf F}(\eta_\varepsilon \uph \delta) = 1 - f_\varepsilon
(\delta)$.

Now define $\eta \in {}^\delta2$ by $\eta(\mu \alpha+ \varepsilon)
= \eta_\varepsilon (\alpha)$ for $\alpha < \lambda, \varepsilon
< \mu.$
\end{enumerate}
Let $E = \{\delta < \lambda : \delta$ is divisible by $\mu^\omega$
and $\varepsilon < \mu \wto \delta \in E_\varepsilon$ and
$(\pa \alpha < \delta) [ \eta \uph \alpha \in \{\eta_i:i <
\delta \}] \}.$
Clearly $E$ is a club of $\lambda$ hence we can find $\delta^*
\in E \cap S$. So by the definition of $\sP_\delta$ we have
$\eta \uph \delta \in \sP_\delta$ and for $\varepsilon <\mu$ we
have $g_{\eta \uph \delta, \varepsilon} \in {}^\delta2$ is equal to
$\eta_\varepsilon \uph \delta$ (Why? note that $\mu \delta=
\mu$ as $\delta \in E$ and see the definition of
$g_{\eta \uph \delta,\varepsilon}$ and of $\eta$, so : $\alpha <
\delta \wto g_{\eta \uph \delta, \varepsilon} (\alpha)
= \eta (\mu \alpha + \varepsilon) = \eta_\varepsilon (\alpha)).$
Hence $h_{\eta \uph \delta} \in {}^\mu2$ is well defined and
by the choice of $\eta$ we have $\varepsilon < \mu \wto
g_{\eta \uph \delta,\varepsilon}=\eta_\varepsilon \uph \delta$ so by
its definition, $h_{\eta \uph \delta}$ for
each $\varepsilon < \mu$ satisfies $h_{\eta \uph
\delta} (\varepsilon) = {\bf F}(g_{\eta\uph \delta,\varepsilon})
={\bf F} =(\eta_\varepsilon \uph \delta).$ Now by the choice of
$f_\varepsilon$ we have ${\bf F} (\eta_\varepsilon \uph \delta) =
f_\varepsilon (\delta) = g^*_\delta (\varepsilon)$ so
$h_{\eta \uph \delta} = g^*_\delta,$ but $h_{\eta \uph \delta}
\in \sP_\delta$ whereas we have chosen $g^*_\delta$ \st\
$g^*_\delta \notin \sP_\delta,$ a contradiction.
\ref{x.1}
\end{proof}

We may consider a generalization.
\begin{definition}
\label{x.3}
\begin{enumerate}
\item We say $\bar{C}$ is a $\lambda-$Wd-parameter
if:
\begin{enumerate}
\item[(a)] $\lambda$ is a regular uncountable,
\item[(b)] $S$ a stationary subsets of $\lambda$,
\item[(c)] $\bar{C} = \lk C_\delta : \delta \in S \rk,
C_\delta \ps \delta$
\end{enumerate}
\begin{enumerate}
\item[(1A)] We say $\bar{C}$ is a $(\lambda,\kappa,\chi)$
-Wd-parameter \underline{if} in addition $(\pa \delta \in S) [\cf
(\delta) = \kappa \wedge |C_\delta| < \chi ].$ We may also say
that $\bar{C}$ is $(S, \kappa,\chi)$-parameter.
\end{enumerate}
\item We say that ${\bf F}$ is a $(\bar{C}, \theta)$
-colouring if:\,$\bar{C}$ is a $\lambda$-Wd-parameter and ${\bf F}$
is a function from ${}^{\lambda>}\lambda$ to $\theta$ \st\ :

if $\delta \in S, \quad  \eta_0,\eta_1 \in
{}^{\delta}\delta$ and $\eta_0
\uph C_\delta = \eta_1 \uph C_\delta$
then ${\bf F} (\eta_0) = {\bf F} (\eta_1)$.
\begin{enumerate}
\item[(2A)] If $\theta=2$ we may omit it writing
$\bar{C}$- colouring
\item[(2B)] In part (2) we can replace ${\bf F}$ by
$\lk F_\delta : \delta \in S \rk$ where $F_\delta :
{}^{(C_\delta)}\delta \to \theta$ \st\  $\eta \in
{}^{\delta}\delta \wedge \delta \in S \to {\bf F} (\eta) =
F_\delta (\eta \uph C_\delta).$ So abusing notation we may write
${\bf F} (\eta \uph C_\delta)$
\end{enumerate}
\item Assume ${\bf F}$ is a $(\bar{C}, \theta)$-clouring,
$\bar{C}$ a $\lambda$-Wd-parameter.

We say $\bar{c} \in {}^S\theta$ (or $\bar{c} \in
{}^{\lambda}\theta$) is an ${\bf F}$-wd-\seq\  if :

(*) for every $\eta \in {}^{\lambda}\lambda$,
the set $\{ \delta \in S : {\bf F} (\eta \uph \delta)=
c_\delta \}$ is a stationary subset of $\lambda$.

We also may say $\bar{C}$ is an $({\bf F},S)$-Wd-\seq.
\begin{enumerate}
\item[(3A)] We say $\bar{c} \in {}^S\theta$ is a
$D-{\bf F}$-Wd-\seq\ if $D$ is a filter on $\lambda$ to which $S$
belongs and

(*)for every $\eta \in {}^\lambda\lambda$ we have
$$
\{ \delta \in S: {\bf F} (\eta \uph \delta) = c_\delta \} \not=
\emptyset \, {\rm mod} D
$$
\end{enumerate}
\item We say $\bar{C}$ is a good $\lambda$-Wd-parameter,
\underline{if} for every $\alpha < \lambda$ we have $\lambda
> |\{C_\delta \cap  \alpha : \delta \in S \}|$.
\end{enumerate}
${}$

Similarly to \ref{x.1} we have
\end{definition}

\begin{claim}
\label{x.4}
Assume
\begin{enumerate}
\item[(a)] $\bar{C}$ is a good $(\lambda,
\kappa, \chi)$-Wd-parameter.
\item[(b)] $|\alpha|^{\lk {\rm tr}, \kappa \rk} < \lambda$
for every $\alpha < \lambda$.
\item[(c)] $\lambda = 2^\mu$ and
$\lambda = \lambda^{< \chi}$
\item[(d)] ${\bf F}$ is a $\bar{C}$- colouring.
\end{enumerate}

\underline{Then} there is a ${\bf F}$-Wd-\seq.
\end{claim}

\begin{proof}
Let ${\rm cd}$ be a 1-to-1 function from
${}^{\mu}\lambda$ onto $\lambda$, for simplicity, and
\Wlog
$$
\alpha = {\rm cd} (\lk \alpha_\varepsilon : \varepsilon
< \mu \rk) \wto \alpha \ge \sup \{ \alpha_\varepsilon
: \varepsilon < \mu \}
$$
and let the function ${\rm cd}_i : \lambda \to \lambda$ for
$i< \mu$ be \st\ ${\rm cd}_i ( \lk {\rm cd}
(\alpha_\varepsilon : \varepsilon < \mu ) \rk ) =
\alpha_i$.

Let $T= \{ \eta:$ for some $C \ps \lambda$
of cardinality $< \chi$, we have $\eta \in
{}^C\lambda \},$ so by assumption (c) clearly
$|T| = \lambda$, so let us list $T$ as
$\{ \eta_\alpha : \alpha < \lambda \}$ with no repetitions,
and let $T_{< \alpha} = \{ \eta_\beta : \beta < \alpha \}$.
For $\delta \in S$ let $\sP_\delta = \{\eta : \eta$ a
function from $C_\delta$ to $\delta$ \st\ for every
$\alpha \in C_\delta$ we have $\eta \uph
(C_\delta \cap \alpha) \in T_{< \delta}.$

By clause (b) of the assumption necessarily $\sP_\delta$
has cardinality $< \lambda$.
For each $\eta \in \sP_\delta$
and $\varepsilon < \mu$ we define $\nu_{\eta,\varepsilon}
\in {}^{C_\delta}\delta$ by $\nu_{\eta,\varepsilon} (\alpha)
= {\rm cd}_\varepsilon (\eta (\alpha))$ for $\alpha \in
C_\delta$. Now for $\eta \in \sP_\delta$, clearly
$\rho_\eta =: \lk {\bf F} (\nu_{\eta, \varepsilon}) :
\varepsilon < \mu\rk$ belongs to ${}^\mu2$. Clearly
$\{ \rho_\eta: \eta \in \sP_\delta\}$ is a subset of
${}^\mu2$ of cardinality $\le |\sP_\delta|$ which as
said above is $< \lambda$. But $|{}^\mu2| = 2^\mu = \lambda$
by clause (c) of the assumption, so we can find $\rho^*_\delta \in
{}^\mu2 \setminus \{ \rho_\eta : \eta \in \sP_\delta\}.$

For each $\varepsilon<\mu$ we can consider the \seq\
$\bar{c}^\varepsilon = \lk 1 - \rho^*_\delta (\varepsilon)
: \delta \in S \rk$ as a candidate for being an ${\bf F}$
-Wd-\seq. If one of then is, we are done. So assume toward
contradiction that for each $\varepsilon< \mu$ there is
$\eta_\varepsilon \in {}^{\lambda}\lambda$ which exemplify
its failure, so there is a club $E_\varepsilon$ of
$\lambda$ \st
\begin{enumerate}
\item[$\boxtimes_1$] $\delta \in S \cap E_\varepsilon
\wto {\bf F} (\eta_\varepsilon \uph C_\delta) \not=
c^\varepsilon_\delta$

and \Wlog
\item[$\boxtimes_2$] $\alpha < \delta \in E_\varepsilon
\wto \eta_\varepsilon (\alpha) < \delta.$

But $c^\varepsilon_\delta = 1 - \rho^*_\delta
(\varepsilon)$ and so $z \in \{0,1\} \, \hbox{\&} \,
z \not= c^\varepsilon_\delta \wto z = \rho_\delta (\varepsilon)$
hence we have gotten
\item[$\boxtimes_3$] $\delta \in S \cap E \wto
{\bf F} (\eta_\varepsilon \uph C_\delta) =
\rho^*_\delta (\varepsilon)$
\end{enumerate}

Define $\eta^* \in {}^{\lambda}\lambda$ by
$\eta^* (\alpha)= {\rm cd} ( \lk \eta_\varepsilon
(\alpha) : \varepsilon < \mu \rk )$, now as
$\lambda$ is regular uncountable clearly
$E=: \{ \delta < \lambda :$ for every
$\alpha < \delta$ we have $\eta^* (\alpha)< \delta$ and if $\delta'
\in S, C' = C_\delta' \cap \alpha$ then $\eta^* \uph C'
\in T_{< \delta} \}$ is a club of $\lambda$
(see the choice of $T, T_{<\delta}$ recall that by assumption (a)
the \seq\ $\bar{C}$ is good, see Definition \ref{x.3}(4)).

Clearly $E^* = \cap \{E_\varepsilon : \varepsilon < \mu \}
\cap E$ is a club of $\lambda$. Now for each $\delta \in
 E^*$, clearly $\eta^* \uph C_\delta \in
\sP_\delta$; just check the definitions of $\sP_\delta$
and $E, E^*$. Now recall $\nu_{\eta^* \uph C_\delta, \varepsilon}$
is the function from $C_\delta$ to $\{0,1\}$ defined by
$$
\nu_{\eta^* \uph C_\delta, \varepsilon} (\alpha)=
{\rm cd}_\varepsilon (\eta^* (\alpha)).
$$
But by our choice of $\eta^*$ clearly ${\rm cd}_\varepsilon
(\alpha))= \eta_\varepsilon (\alpha)$, so
$$
\alpha \in C_\delta \wto \nu_{\eta^* \uph C_\delta, \varepsilon}
(\alpha) = \eta_\varepsilon (\alpha) \quad \hbox{so} \quad
\nu_{\eta^* \uph C_\delta, \varepsilon} =
\eta_\varepsilon \uph C_\delta ,
$$
Hence ${\bf F} (\nu_{\eta^* \uph C_{\delta,\varepsilon}})
= {\bf F} (\eta_\varepsilon \uph C_\delta)$, however as
$\delta \in E^* \ps E_\varepsilon$ clearly ${\bf F}
(\eta_\varepsilon \uph C_\delta) = \rho^*_\delta
(\varepsilon)$, together ${\bf F} (\nu_{\eta^* \uph
C_\delta, \varepsilon} ) = \rho^*_\delta (\varepsilon)$.

As $\eta^* \uph C_\delta \in \sP_\delta$ clearly $\rho_{\eta^* \uph
c_\delta} \in {}^{\mu}2$, moreover for each $\varepsilon < \mu$
 we have  $\rho_{\eta^* \uph C_\delta} (\varepsilon)$, see its
definition above, is equal to ${\bf F} (\nu_{\eta^* \uph
C_\delta, \varepsilon})$ which by the previous sentence
is equal to $\rho^*_\delta (\varepsilon)$. As this holds for every
$\varepsilon < \mu$ and $ \rho_{\eta^* \uph C_\delta},
\rho^*_\delta$ are members of ${}^\mu2$, clearly they are
equal. But $\eta^* \uph C_\delta \in \sP_\delta$ so
$\rho_{\eta^* \uph C_\delta} \in \{ \rho_\eta : \eta
\in \sP_\delta\}$ whereas $\rho^*_\delta$ has been chosen
outside this set, contradiction.
\end{proof}

Well, are there good $(\lambda, \kappa, \kappa)$-parameter?
(on $I[\lambda]$ see \cite[\S 1]{Sh:420}).

\begin{claim}
\label{x.5}
\begin{enumerate}
\item If $S$ is a stationary subset of the regular cardinal
$\lambda$ and $S \in I [\lambda]$ and $(\pa \delta \in S)
\cf (\delta) = \kappa$ \underline{then} for some club $E$ of
$\lambda$, there is a  good $(S \cap E, \kappa, \kappa)$-parameter.
\item If $\kappa = \cf(\kappa), \kappa^+ < \lambda=\cf (\lambda)$
\underline{then} there is a  stationary $S\in I [\lambda]$ with $(\pa
\delta \in S) [\cf (\delta) = \kappa]$.
\end{enumerate}
\end{claim}

\begin{proof}
\begin{enumerate}
\item By the definition of $I [\lambda]$
\item By \cite[\S 1]{Sh:420}.
\end{enumerate}
\end{proof}

We can note
\begin{claim}
\label{x.6}
\begin{enumerate}
\item Assume the assumption of \ref{x.4} or
\ref{x.1} with $C_\delta = \delta$ and $D$
is a $\mu^+$- complete filter on $\lambda, S \in D$,
and $D$ include the club filter. \underline{Then} we can get that
there is a $D-{\bf F}$-Wd-\seq.
\item In \ref{x.4}, we can weaken the demand
$\lambda = 2^\mu$ to $\lambda = \cf (2^\mu)$
that is, assume
\begin{enumerate}
\item[(a)] $\bar{C}$ is a good ($\lambda, \kappa,
\chi)$-Wd-parameter.
\item[(b)] $|\alpha|^{\lk{\rm tr},\kappa \rk} < 2^\mu$ for every
$\alpha<\lambda$.
\item[(c)]$\lambda=\cf(2^\mu)$ and $2^\mu = (2^\mu)^{<\chi}$
\item[(d)] ${\bf F}$ is a $\bar{C}$-colouring
\item[(e)]$D$ is a $\mu^+$-complete filter on $\lambda$ extending the
club filler to which ${\rm Dom}(\bar{C})$ belongs.

\underline{Then} there is a $D-{\bf F}$-Wd-\seq.
\end{enumerate}
\item In \ref{x.4}+\ref{x.6}(2) we can omit ``$\lambda$ regular''.
\end{enumerate}
\end{claim}

\begin{proof}
\begin{enumerate}
\item The same proof.
\item Let $H^* : \lambda \to 2^\mu$ be \incr\  \cont\  with
unbounded range and let $S \in I[\lambda]$ be stationary, \st\
$(\pa \delta \in S) \cf (\delta)=\kappa,$ and $\bar{C}= \lk C_\delta:
\delta \in S \rk$ is a good $(\cf (\lambda), \kappa, \kappa)$-Wd-
parameter,let
$$
S' = \{h^* (\alpha): \alpha \in S\},\, C'_{h^*(\delta)}=
\{h^* (\alpha) : \alpha \in C_\delta \}, \bar{C}'=\lk C_\beta
:\beta \in S' \rk
$$
and repeat the proof using $\lambda' = 2^\mu, \bar{C}'=
\lk C'_\delta : \delta \in S' \rk$ instead $\lambda, \bar{C}$.
Except that in the choice of the club $E$ we should use
$E' = \{ \delta < \lambda$: for every $\alpha \in \delta \cap$
Rang $(h^*)$ we have $\eta^* (\alpha) < \delta$ and
$\delta$ is a limit ordinal and $ \delta' \in S' \wedge
C'= C'_\delta \cap \alpha \wto \eta^* \uph C'
\in T_{< \delta}\}$.
\item Similarly.
\end{enumerate}
\end{proof}

This lead to considering the natural related ideal.

\begin{definition}
\label{x.7}
Let $\bar{C}$ be a $(\lambda, \kappa, \chi)$-
parameter.
\begin{enumerate}
\item For a family $\sF$ of $\bar{C}$-colouring and
$\sP \ps {}^\lambda 2,$ let
${\rm id}_{\bar{C}, \sF,\sP}$ be

$\{W \ps \lambda:$ for some
${\bf F} \in \sF$ for every $\bar{c} \in \sP$ for some
$\eta \in {}^\lambda\lambda$ the set
$$
\{ \delta \in W \cap S : {\bf F} (\eta \uph C_\delta) = c_\delta\}
\quad \hbox{is not stationary}\}.
$$
\item If $\sP$ is the family of all $\bar{C}$-
colouring we may omit it.
If we write Def instead $\sF$ this mean as in
\cite[\S 1]{Sh:576}.
\end{enumerate}

${}$

We can strengthen \ref{x.4} as follows.
\end{definition}

\begin{definition}
\label{x.8}
We say the $\lambda$-colouring ${\bf F}$ is
$(S, \chi)$- good if:
\begin{enumerate}
\item[(a)] $S \ps \{ \delta < \lambda : \cf (\delta)
< \chi \}$ is stationary
\item[(b)] we can find $E$ and $\lk C_\delta : \delta
\in S \cap E \rk$ \st\
\begin{enumerate}
\item[$(\alpha)$] $E$ a club of $\lambda$.
\item[$(\beta)$] $C_\delta$ is an unbounded subset of
$\delta, |C_\delta| < \chi$.
\item[$(\gamma)$] if
$ \rho,\rho' \in {}^{\delta}\delta,\, \delta \in S \cap E,$
\quad and \quad
$\rho' \uph C_\delta = f \uph C_\delta$

then \quad ${\bf F} (\rho') = {\bf F} (\rho)$
\item[$(\delta)$] for every $\alpha < \lambda$ we have
$$
 \lambda > | \{C_\delta \cap
\alpha: \delta \in S \cap E \}|
$$
\item[$(\epsilon)$] $\delta \in S \wto |\delta|^{\lk {\rm tr},\cf
(\delta)\rk} <\lambda$ or just $\delta \in S \wto \lambda > |\{
C:C \ps \delta$ is unbound and for every $\alpha< \delta$ for some
$\gamma \in S$ we have $C \cap \alpha = C_\gamma \cap \alpha]:$
\end{enumerate}
\end{enumerate}
\end{definition}

\begin{claim}
\label{x.9}
Assume
\begin{enumerate}
\item[(a)] $\lambda = \cf(2^\mu)$
\item[(b)] ${\bf F}$ is an $(S, \kappa)$- good
$\lambda$-colouring.
\end{enumerate}

\underline{Then} there is a $({\bf F},S)$-Wd-\seq\  see Definition
\ref{x.3}(3).
\end{claim}

\begin{remark}
\label{x.10}
So if $\lambda=\cf(2^\mu)$ and we let $\Theta_\lambda =:\{
\theta=\cf(\theta)$ and $(\pa \alpha<\lambda) (|\alpha|^{\lk {\rm tr},
\theta\rk} < \lambda)\}$ then
\begin{enumerate}
\item[(a)] $\Theta_\lambda$ ``large'' (e.g. contains every large
enough $\theta \in {\rm Reg} \cap \beth_\omega$ if
$\beth_\omega<\lambda)$ and
\item[(b)] if $\theta \in \Theta_\lambda \wedge \theta^+ < \lambda$
then there is a stationary $S \in I[\lambda]$ \st\ $\delta \in S \wto
\cf (\delta)=\theta.$
\item[(c)] if $\theta,S$ are as above \underline{then}
there is a good $\lk C_\delta:\delta\in S\rk$
\item[(d)] for $\theta,S,\bar{C}$ as above, if ${\bf F}
= \lk F_\delta:\delta \in S\rk$ and $F_\delta(\eta)$ depend just
on $\eta \uph C_\delta$ and $D$ is a normal ultrafilter on $\lambda$
(or less), and lastly $S \in D$ \underline{then} there is an
$D-{\bf F}$-Wd-\seq; see Definition \ref{x.3}(3A).
\end{enumerate}
\end{remark}

\section{On versions of precipitousness}

\begin{definition}
\label{2.1}
\begin{enumerate}
\item We say the $D$ is $(\bbP,\name{D})$
-precipituous if
\begin{enumerate}
\item[(a)] $D$ is a normal filter on $\lambda$, a regular
uncountable cardinal.
\item[(b)] $\bbP$ is forcing notion with $\emptyset_\bbP$
minimal.
\item[(c)] $\name{D}$ a $\bbP$-name of an
ultrafilter of the Boolean Algebra \, $\sP(\lambda)$
\item[(d)] letting for $p \in \bbP$
$$
D_{p,\name{D}} =: \{A \ps \lambda: p \Vdash
A \in \name{D} \}
$$
we have:
\begin{enumerate}
\item[$(\alpha)$] $D_{\emptyset_\bbP, \name{D}}=
D$ and
\item[$(\beta)$] $D_{p,\name{D}}$ is normal filter
on $\lambda$
\end{enumerate}
\item[(e)] $\Vdash_\bbP \quad \hbox{``} {\bf V}^\lambda /
\name{D}$ is well founded''.
\end{enumerate}
(1A) If $\name{D}$ is clear from the context (as in
part (2)) we may  omit $\name{D}$.
\item For $\lambda$ regular uncountable and $D$ a normal filter on
$\lambda$ let ${\rm NOR}_D = \{D' : D'$ a normal filter on
$\lambda$ extending $D \}$ ordered by inclusion and
$\name{D} = \cup \{D' : D' \in
\name{G}_{{\rm NOR}_D} \}$
\end{enumerate}
Woodin [W99] define and was be
interested in semi saturation for $\lambda = \aleph_2$.
\end{definition}

\begin{definition}
\label{2.2}
For $\lambda$ regular uncountable cardinal,
a normal filter $D$ on $\lambda$   is called semi saturated when for
every forcing notion $\bbP$ and $\bbP$-name
$\name{D}$ of a normal (for regressive $f \in {\bf V})$
ultrafilter on $\sP (\lambda)^{\bf V}$, we have:
$D$ is $(\bbP, \name{D})$- precipitous.

Woodin prove ${\rm Con} (D_{\omega_2} \uph S^2_0$ is semi
saturated). He proved that the existence of such filter has
large consistency strength by proving 2.3 below.
 This is related to \cite[V]{Sh:g}.
\end{definition}

\begin{claim}
\label{2.3}
If $\lambda = \mu^+,\, D$ a semi saturated filter or
$\lambda$, \underline{then} every $f \in {}^{\lambda}\lambda$ is
$<_D$- than the $\alpha$-th function for some $\alpha<\lambda^+$
(on the $\alpha$-th function see e.g \cite[XVII, \S 3]{Sh:g})

${}$

In fact
\end{claim}

\begin{claim}
\label{2.5x}
If $\lambda = \mu^+,\, D$ is ${\rm NOR}_\lambda$-precipitous
\underline{then} every $f \in {}^{\lambda}\lambda$ is
$<_D$- smaller than the $\alpha$-th function for some
$\alpha< \lambda^+$
\end{claim}

\begin{proof}
The point is that
\begin{enumerate}
\item[(a)] if $D$ is a normal filter on $\lambda, \lk f_\alpha:\alpha
< \lambda^+ \rk$ is ${}^<D$ -\incr\ in $\lambda$ and $f \in
{}^\lambda\lambda, \alpha<\lambda^+ \wto \neg (f \le_D f_\alpha)$
\underline{then} there is a normal filter $D^+_1$ on $\lambda$
extending $D$ \st\ $\alpha<\lambda^+ \wto f_\alpha <_{D_1} f$
\item[(b)] if $\lk f_\alpha:\alpha \le \lambda^+ \rk$ is $<_D$-
\incr\ $f_\alpha \in {}^\lambda\lambda$, and $\lambda=\mu^+$ and
$X=\{\delta<\lambda:\cf (f_{\lambda^+} (\delta))= \theta\}
\not= \emptyset\, {\rm mod} D$  then there are functions
$g_i \in {}^\lambda\lambda$ for $i<\theta$
\st\ $g_i< f_{\lambda^+} {\rm mod} \,(D+X)$, and
$(\pa \alpha < \lambda^+) (\exists i < \theta)
(\neg g_i <_D f_\alpha).$
\end{enumerate}
\end{proof}

\begin{claim}
\label{2.5}
\begin{enumerate}
\item If $\lambda= \mu^+ \ge \beth_\omega$ \underline{then}
 the club filter on
$\lambda$ is not semi saturated.
\item If $\lambda= \mu^+ \ge \beth_\omega$ \underline{then}
for every large enough regular $\kappa< \beth_\omega$, there is
no semi saturated normal filter $D^*$ on $\lambda$ to which
$S^\lambda_\kappa = \{\delta < \lambda : \cf (\delta)=\kappa \}$
belongs.
\item If $\lambda=\mu^+ > \kappa = \cf (\kappa) > \aleph_0$
and for every $f \in {}^\kappa\lambda$ we have
${\rm rk}_{J^{bd}_\kappa} (f) < \lambda$
\underline{then} there is no semi saturated normal filter
$D^*$ on $\lambda$ to which $\{\delta < \lambda : \cf
(\delta)=\kappa \}$ belongs.
\item In 1), 2), 3), if ``$D$ is ${\rm Nor}_D$-semi saturated''
\underline{then} the conclusion holds for $D$.
\end{enumerate}
\end{claim}

\noindent{\sc Remark:}\quad We can replace $\beth_\omega$ by
any strong limit uncountable cardinal.
\medskip

\begin{proof}
\begin{enumerate}
\item Follows by (2)
\item By \cite{Sh:460} for some $\kappa_0 < \beth_\omega$,
for every regular $\kappa \in (\kappa_0, \beth_\omega)$
we have: $\mu^{\lk \kappa \rk} = \mu$. Let $D=\{A \ps \kappa
: \sup (\kappa \setminus A)<\kappa\}.$

By part (3) it is enough to prove

$\boxtimes$\quad if $f \in {}^\kappa\lambda$ then ${\rm rk}_D
(f) < \lambda$

\underline{proof of $\boxtimes$} If not then for every $\alpha<
\lambda$ there is
$$
f_\alpha \in {}^\kappa\lambda \quad \hbox{\st}\quad f_\alpha <_D f \quad
\hbox{and}\quad {\rm rk}_D (f) = \alpha
$$
and define
$$
D_\alpha = : \{A \ps \kappa: \,A \in D\quad \hbox{\underline{or}}
\quad \kappa \setminus A \notin D, \hbox{and} \quad {\rm rk}_{D+
(\kappa\setminus A)} (f_\alpha) < \alpha \}.
$$
This is a $\kappa$-complete filter on $\kappa$ see \cite{Sh:589}.
So for some $D^*$ the set $A=\{\alpha:D_\alpha =D^*\}$
is unbounded in $\lambda$.
By \cite[\S 4x]{Sh:589} (alternatively use \cite[V]{Sh:g} on normal
filters)

(*) for $\alpha<\beta$ from $A, f_\alpha <_{D^*} f_\beta$
and $D^*$ is a $\kappa$-complete filter on $\kappa$.

But as $\mu=\mu^{\lk \kappa \rk}$ letting $\alpha^*=\sup ({\rm Rang}
(f)) +1$ which is $<\lambda$, so $|\alpha^*| \le \mu$,
there is a family $\sP \ps [\alpha^*]^\kappa$ \st\ for every
$a \in [\alpha^*]^\kappa$, for some $i (*) < \kappa$ and
$a_i \in \sP$ for $i<i(*)$ we have $a \ps \bigcup\limits_{i<i(*)} a_i$
hence for every $\alpha \in A$, for some $a_\alpha \in \sP$ we have
$$
 \{i<\kappa:f_\alpha(i) \in a_\alpha \}
\not= \emptyset \, {\rm mod} \, D^*.
$$
So for some $a^*$ and unbounded $B \ps A$ we have $\alpha \in
B \wto  a_\alpha = a^*$ and moreover for some $b^* \ps \kappa$ we have
$\alpha \in B \wto b^* = \{i < \kappa: f_\alpha (i) \in a^* \}$
and moreover $\alpha \in B \wto f_\alpha \uph b^* = f^*$.
But this contradict (*).
\item We can find $\lk u_{\alpha,\varepsilon}:
\varepsilon < \lambda, \alpha < \lambda^+ \rk$ \st :
\begin{enumerate}
\item[(a)] $\lk u_\alpha,\varepsilon:\varepsilon < \lambda \rk$ is
$\ps$-increasing \cont\ \st\ $|u_{\alpha,\varepsilon}| < \lambda,$
and $\cup \{u_{\alpha,\varepsilon}:\varepsilon < \lambda \} = \alpha$.
\item[(b)] if $\alpha<\beta<\lambda^+$ and $\alpha \in
u_{\beta,\varepsilon}$ then $u_{\beta,\varepsilon}
\cap \alpha= u_{\alpha,\varepsilon}$.
\end{enumerate}

Let $f_\alpha \in {}^\lambda\lambda$ be $f_\alpha (\varepsilon)=
{\rm otp} (u_{\alpha,\varepsilon})$, so it is well known that
$f_\alpha / D_\lambda$ is the $\alpha$-th function, in particular
$\alpha<\beta \wto f_\alpha <_{D_\lambda} f_\beta$ where
$D_\lambda$ is the club filter on $\lambda$; in fact
$\alpha<\beta<\lambda^+ \wto f_\alpha<_{J^{bd}_\lambda}
f_\beta.$
Choose\footnote{recall $S^\lambda_\kappa=\{\delta<\lambda:
\cf(\delta)_\kappa\}$}
 $\bar{C} = \lk C_\delta:\delta \in S^\lambda_\kappa
\rk,\,  C_\delta$ a club of $\delta$ of order type $\kappa$, and let
$g_\delta \in {}^\kappa\delta$ enumerate $C_\delta$, i.e.
$g_\delta (i)$ is the $i$-th member of $C_\delta$

For $\zeta<\lambda$ let $g^*_\zeta \in {}^\kappa\lambda$ be
constantly $\zeta$, and let
$g^* \in {}^\lambda\lambda$ be defined by $g^*(\zeta)= {\rm
rk}_{J^{bd}_\kappa} (g^*_\zeta)$
\begin{enumerate}
\item[$(*)_0$] $g^* \in {}^\lambda\lambda$

[why? by an assumption]

For $\alpha < \lambda^+$ we define $f^*_\alpha \in
{}^\lambda\lambda$ by:
\begin{equation*}
f^*_\alpha (\varepsilon) =
\begin{cases} {\rm rk}_{J^{bd}_\kappa} (f_\alpha \circ
g_\epsilon)\quad \hbox{if} \quad \varepsilon \in S^\lambda_\kappa
\\ 0 \quad \quad\quad \quad \quad \quad \,
\hbox{if} \quad  \varepsilon \in \lambda
\setminus S^\lambda_\kappa
\end{cases}
\end{equation*}
Note that $f_\alpha \circ g_\delta$ is a function from $\kappa$
to $\lambda$.

Now
\item[$(*)_1$] $f^*_\alpha \in {}^\lambda\lambda$ for
$\alpha < \lambda^+$

[Why? as $f_\alpha \circ g_\delta \in {}^\kappa\lambda$, so by a
hypothesis ${\rm rk}_{J^{bd}_\kappa} (f_\alpha \circ g_\delta)
< \lambda ]$
\item[$(*)_2$] for $\alpha < \lambda^*$
$$
(*)_\alpha^2 \, E_\alpha =\{\delta < \lambda:\,
\hbox{if}\quad \varepsilon < \lambda \quad \hbox{then} \quad
f^*_\alpha (\varepsilon) < \delta \}
$$
is a club of $\lambda$

[Why? Obvious]
\item[$(*)_3$] for $\alpha < \lambda^+$ we have
$$
\delta \in E_\alpha \wto f^*_\alpha (\delta) <
g^* (\delta), \, \hbox{so} \quad f^*_\alpha <_{D_\lambda}
g^* \in {}^\lambda\lambda
$$
[Why? the first statement by the definition of $E_\alpha$ and
of $g^* (\delta)$. The second by the first $(*)_0$.]
\item[$(*)_4$] if $\alpha < \beta < \lambda^+$
\underline{then} $f^*_\alpha <_{J^{bd}_\lambda} f^*_\beta$
hence $f^*_\alpha <_{D_\lambda} f^*_\beta$

[Why? the first as $f_\alpha <_{J^{bd}_\lambda} f_\beta$
hence for some $\varepsilon < \lambda$, we have
\begin{align*}
\varepsilon < \zeta <\lambda \to f_\alpha (\zeta) < f_\beta
(\zeta) \,\hbox{hence} \quad \delta \in S^\lambda_\kappa \setminus
(\varepsilon+1) \wto  \\
f_\alpha \uph C_\delta
<_{J^{bd}_{C_\delta}} f_\beta \uph C_\delta \wto
f_\alpha \circ g_\delta <_{J^{bd}_\kappa} f_\beta \circ
g_\delta \wto
\end{align*}
\end{enumerate}
Let $f^*_{\lambda^+} = : g^*$, so

$$
(*) \, \alpha \le \lambda^+ \wto f^*_\alpha \in {}^\lambda\lambda
\quad \hbox{and} \quad \alpha < \beta \le \lambda^+
\wto f_\alpha <_{D_\lambda} f_\beta
$$
This of course suffice by \ref{2.3}.
\end{enumerate}
\end{proof}

\noindent{\sc Remark:}
\label{2.5a}
In the proof of \ref{2.5}(2) it is enough that ${\bf U}_{J^{bd}_\kappa}
(\mu) = \mu$ (see \cite{Sh:589}).
\medskip


\par\noindent [References of the form {\tt math.XX/$\cdots$} refer to the
{\tt xxx.lanl.gov} archive]  \par

\end{document}